\documentclass[twoside,a4paper]{article}

\usepackage{amsmath,amsthm,amssymb,latexsym}

\theoremstyle{plain}%
\newtheorem{proposition}{Proposition}[section]%
\newtheorem{theorem}[proposition]{Theorem}%
\newtheorem{corollary}[proposition]{Corollary}%
\theoremstyle{definition}%
\newtheorem{definition}[proposition]{Definition}%

\newcommand{\wap}{\operatorname{WAP}}
\newcommand{\ap}{\operatorname{AP}}

\newcommand{\FIN}{\operatorname{FIN}}

\newcommand{\proten}{{\widehat{\otimes}}}
\newcommand{\mc}[1]{\mathcal{#1}}
\newcommand{\ip}[2]{{\langle {#1} , {#2} \rangle}}

\newcommand{\aone}{\Box}
\newcommand{\atwo}{\Diamond}
\newcommand{\newprod}{\underset{K,\mc U}{\star}}
\newcommand{\newprodinline}{\text{\raisebox{.7ex}{$\newprod$}}}

\begin{document}

\large
\title{Ultrapowers of Banach algebras and modules}
\author{Matthew Daws}
\maketitle

\begin{abstract}
The Arens products are the standard way of extending the product
from a Banach algebra $\mc A$ to its bidual $\mc A''$.  Ultrapowers
provide another method which is more symmetric, but one that in general
will only give a bilinear map, which may not be associative.  We show
that if $\mc A$ is Arens regular, then there is at least one way to
use an ultrapower to recover the Arens product, a result previously
known for C$^*$-algebras.  Our main tool is a Principle of Local
Reflexivity result for modules and algebras.

\medskip


\end{abstract}

\section{Introduction}

Given a Banach algebra $\mc A$, there are two canonical ways
to extend the products on $\mc A$ to products on the bidual
$\mc A''$, called the Arens products (see \cite{Arens}).
Another, less common, way to view the bidual of a Banach space
is as a quotient of an ultrapower of that space (see \cite{Hein}).
As the ultrapower of a Banach algebra is again a Banach algebra,
this suggests another way of defining a ``product'' on $\mc A''$.

We recall below the notion of an ultrapower of $\mc A$ with respect
to an ultrafilter $\mc U$, written $(\mc A)_{\mc U}$.  By
weak$^*$-compactness of the ball of $\mc A''$, we can always define
a norm-decreasing map
\[ \sigma_{\mc U}:(\mc A)_{\mc U}\rightarrow\mc A'';\quad
\ip{\sigma_{\mc U}((a_i))}{\mu} = \lim_{i\rightarrow\mc U} \ip{\mu}{a_i}
\qquad (\mu\in\mc A', (a_i)\in(\mc A)_{\mc U}). \]
It follows from the Principle of Local Reflexivity
(see \cite[Proposition~6.7]{Hein} or Theorem~\ref{PLR} below) there for
a suitable $\mc U$, there is an isometry $K:\mc A''\rightarrow
(\mc A)_{\mc U}$ such that $\sigma_{\mc U}\circ K$ is the identity on
$\mc A''$.  Using this, we may define a bilinear map
\[ \star = \newprod : \mc A''\times\mc A'' \rightarrow \mc A'';
\quad \Phi\star\Psi = \sigma_{\mc U}( K(\Phi) K(\Psi) ) \qquad
(\Phi,\Psi\in\mc A''). \]
This idea was explored by Godefroy, Iochum and Loupias in
\cite{GI} and \cite{IL}. In general this might only lead to a
bilinear map (which can fail to be associative).  However, for
C$^*$-algebras, we always recover the Arens products.

In this paper, we shall explore these ideas further, and also
consider related ideas for Banach modules.  The Banach space
tool which relates an ultrapower of a Banach space $E$ to
the bidual $E''$ is the Principle of Local Reflexivity.  It
should hence come as no surprise that our line of attack is to
prove various strengthenings of the Principle of Local Reflexivity
for Banach algebras and modules.  For example, we show that when
$\mc A$ is Arens regular, that is, the two Arens products agree
on $\mc A''$, then there is at least one way to use an ultrapower
of $\mc A$ to induce the same product on $\mc A''$.  As a result,
we get a ``symmetric'' definition of the product on $\mc A''$.

We introduce some notation below, and summarise the results of
\cite{GI} and \cite{IL}.  We then study the
Principle of Local Reflexivity, and prove versions for Banach
modules and Banach algebras which allow us to draw conclusions
about ultrapowers.  We investigate when $\star$ can actually
be chosen to be an algebra homomorphism, and how these ideas
interact with dual Banach algebras.  Finally, we show that
$\star$ can be badly behaved even for C$^*$-algebras if the
map $K$ is not chosen to be an isometry.

\subsection{Notation and basic concepts}\label{Intro}

Let $E$ be a Banach space.  We write $E'$ for the dual space of $E$,
and for $x\in E$ and $\mu\in E'$, we write $\ip{\mu}{x}$ for $\mu(x)$.
Recall the canonical map $\kappa_E:E\rightarrow E''$ defined by
$\ip{\kappa_E(x)}{\mu} = \ip{\mu}{x}$ for $x\in E$ and $\mu\in E'$.
When $\kappa_E$ is an isomorphism, we say that $E$ is \emph{reflexive}.

Recall the notions of filter and ultrafilter.  Let $\mc U$ be a
non-principal ultrafilter on a set $I$, and let $E$ be a Banach
space.  We form the Banach space
\[ \ell^\infty(E,I) = \Big\{ (x_i)_{i\in I}\subseteq E : \|(x_i)\| :=
\sup_{i\in I} \|x_i\| < \infty \Big\}, \]
and define the closed subspace
\[ \mc N_{\mc U} = \Big\{ (x_i)_{i\in I} \in \ell^\infty(E,I) :
\lim_{i\rightarrow\mc U} \|x_i\| =0 \Big\}. \]
Thus we can form the quotient space, called
the \emph{ultrapower of $E$ with respect to $\mc U$},
\[ (E)_{\mc U} := \ell^\infty(E,I) / \mc N_{\mc U}. \]
In general, this space will depend on $\mc U$ (and upon, for example,
if the continuum hypothesis holds), though many properties
of $(E)_{\mc U}$ turn out to be independent of $\mc U$, as long as
$\mc U$ is sufficiently ``large'' in some sense.

We can verify that, if $(x_i)_{i \in I}$ represents an equivalence
class in $(E)_{\mc U}$, then
\[ \| (x_i)_{i\in I} + \mc N_{\mc U} \|
= \lim_{i\rightarrow\mc U} \|x_i\|. \]
We shall abuse notation and write $(x_i)$ for the equivalence
class it represents; of course, it can be checked that any
definition we make is independent of the choice of representative
of equivalence class.  There is a canonical isometry
$E \rightarrow (E)_{\mc U}$ given by sending $x\in E$ to
the constant family $(x)$.  We again abuse notation and write
$x \in (E)_{\mc U}$, identifying $E$ with a closed subspace of $(E)_{\mc U}$.

Recall the notion of a \emph{countably-incomplete} ultrafilter,
for which see \cite{Hein}.  To avoid set-theoretic complications,
we shall henceforth assume that all our ultrafilters are countably-incomplete.

There is a canonical map $(E')_{\mc U} \rightarrow (E)_{\mc U}'$ given by
\[ \ip{(\mu_i)}{(x_i)} = \lim_{i\rightarrow\mc U} \ip{\mu_i}{x_i}
\qquad ( (\mu_i)\in (E')_{\mc U}, (x_i)\in (E)_{\mc U} ). \]
This map is an isometry, and so we identify $(E')_{\mc U}$ with
a closed subspace of $(E)'_{\mc U}$.  It is shown in
\cite[Proposition~7.1]{Hein} that $(E)'_{\mc U} = (E')_{\mc U}$ if and only if
$(E)_{\mc U}$ is reflexive.

For Banach spaces $E$ and $F$, we write $\mc B(E,F)$
for the space of bounded linear operators from $E$ to $F$.  Then
there is a canonical isometric map $( \mc B(E,F) )_{\mc U}
\hookrightarrow \mc B( (E)_{\mc U}, (F)_{\mc U} )$ given by
\[ T(x) = (T_i(x_i)) \qquad ( T=(T_i)\in ( \mc B(E,F) )_{\mc U},
x=(x_i)\in (E)_{\mc U} ). \]
We shall often identify $(\mc B(E,F))_{\mc U}$ with its image
in $\mc B( (E)_{\mc U}, (F)_{\mc U} )$.

When $\mc A$ is a Banach algebra, $(\mc A)_{\mc U}$ becomes
a Banach algebra under the pointwise product.  This follows, as it
is easy to show that $\mc N_{\mc U}$ is a closed ideal in the Banach
algebra $\ell^\infty(\mc A,I)$.  Then $\mc A$ is
commutative if and only if $(\mc A)_{\mc U}$ is; $(\mc A)_{\mc U}$
is unital if $\mc A$ is unital.  If $\mc A$ is a Banach $*$-algebra
(see \cite[Chapter~3]{Dales}) or a C$^*$-algebra, then $(\mc A)_{\mc U}$
is a Banach $*$-algebra or a C$^*$-algebra, respectively, with the
involution defined pointwise.  Thus, as the class of $C(K)$ spaces
for compact, Hausdorff spaces $K$ is the class of commutative,
unital C$^*$-algebras, we see that the ultrapower of a $C(K)$
space is again a $C(K)$ space.

We now recall the Arens products on $\mc A''$.  Firstly, we turn
$\mc A'$ into a $\mc A$-bimodule in the usual fashion,
\[ \ip{a\cdot\mu}{b} = \ip{\mu}{ba}, \quad
\ip{\mu\cdot a}{b} = \ip{\mu}{ab} \qquad (a,b\in\mc A,\mu\in\mc A'). \]
In a similar way, $\mc A''$, $\mc A'''$, and so forth also become
$\mc A$-bimodules.  Then we define bilinear maps
$\mc A''\times\mc A', \mc A'\times\mc A''\rightarrow\mc A'$ by
\[ \ip{\Phi\cdot\mu}{a} = \ip{\Phi}{\mu\cdot a}, \quad
\ip{\mu\cdot\Phi}{a} = \ip{\Phi}{a\cdot\mu} \qquad (\Phi\in\mc A'',
\mu\in\mc A',a\in\mc A). \]
Finally, we define bilinear map $\aone,\atwo:\mc A''\times\mc A''\rightarrow
\mc A''$ by
\[ \ip{\Phi\aone\Psi}{\mu} = \ip{\Phi}{\Psi\cdot\mu},
\quad \ip{\Phi\atwo\Psi}{\mu} = \ip{\Psi}{\mu\cdot\Phi}
\qquad (\Phi,\Psi\in\mc A'', \mu\in\mc A'). \]
These are associative products which extend the natural action
of $\mc A$ on $\mc A''$, called the \emph{first} and \emph{second Arens products}.
See \cite[Section~3.3]{Dales} or \cite[Section~1.4]{Palmer1} for further details.
Thus $\aone$ and $\atwo$ agree with the usual product on $\kappa_{\mc A}(\mc A)$.
When $\aone$ and $\atwo$ agree on all of $\mc A''$, we say that
$\mc A$ is \emph{Arens regular} (see also Section~\ref{wap_ap_sec} below).
As stated above, for suitable $\mc U$, given $\Phi,\Psi\in\mc A''$,
we can find bounded families $(a_i)$ and $(b_i)$ with $(a_i)$ tending
to $\Phi$ weak$^*$ along $\mc U$, and $(b_i)$ tending to $\Psi$.
Then
\[ \ip{\Phi\aone\Psi}{\mu} =
\lim_{j\rightarrow\mc U} \lim_{i\rightarrow\mc U} \ip{\mu}{a_ib_j}, \quad
\ip{\Phi\atwo\Psi}{\mu} =
\lim_{i\rightarrow\mc U} \lim_{j\rightarrow\mc U} \ip{\mu}{a_ib_j}
\qquad (\mu\in\mc A'). \]
We show below that when $\mc A$ is Arens regular, we can find
a more ``symmetric'' version of these formulae.

Recall the map $\star = \newprodinline$ defined above.
We shall henceforth always assume that $K:\mc A''\rightarrow (\mc A)_{\mc U}$
is such that $\sigma\circ K$ is the identity on $\mc A''$ and that
$K\circ\kappa_{\mc A}$ is the canonical map $\mc A\rightarrow
(\mc A)_{\mc U}$.  This is enough to ensure that:
\begin{itemize}
\item For $\Phi\in\mc A''$ and $a\in\mc A$, we have $\Phi\star\kappa_{\mc A}(a) =
\Phi\cdot a$ and $\kappa_{\mc A}(a)\star\Phi = a\cdot\Phi$;
\item if $\mc A$ has a unit $e_{\mc A}$, then $\kappa_{\mc A}(e_{\mc A})$
is a unit for $\star$.
\end{itemize}

In \cite{IL}, the authors make the further assumption that $K$ is
always an \emph{isometry}.  Under this extra condition, from the proof
of \cite[Corollary~II.2]{GI}, it follows that when $\mc A$ is a C$^*$-algebra
the map $\newprodinline$ always agrees with $\aone=\atwo$
(recall that a C$^*$-algebra is always Arens regular, see for
example \cite[Corollary~I.2]{GI}).  We shall show that without
this isometric condition on $K$, we do not always have that
$\newprodinline = \aone$, even for some commutative C$^*$-algebras.
In \cite[Example~3, Page~55]{GI}, it is shown that when $\mc A=\ell^1$
with pointwise product (which is easily seen to be Arens regular),
then for some $K$, we do not have that $\newprodinline$ agrees with the
Arens products.

In \cite[Definition~5]{IL}, the authors say that $\newprodinline$
is \emph{regular} if it is separately weak$^*$-continuous.  They
show in \cite[Proposition~6]{IL} that if $\newprodinline$
is regular for some $K$ and $\mc U$, then $\mc A$ is Arens regular,
and that $\newprodinline=\aone=\atwo$.  This follows fairly
easily from weak$^*$-continuity.  Conversely, in \cite[Theorem~12]{IL},
the authors show that if $\mc A$ is commutative and not Arens regular,
then $\newprodinline$, as it is always a commutative bilinear
map, never agrees with either Arens product.  Let $E$ be a reflexive
Banach space with the \emph{approximation property} (see
\cite[Section~4]{Ryan} or \cite[Section~VIII]{DU}).
Let $\mc A = \mc K(E)$, the algebra of compact operators on $E$.
Then \cite[Corollary~14]{IL} shows that $\newprodinline$
is associative on $\mc A$ if and only if it is regular.  The second
remark after this result in \cite{IL} asks if there is in general
any link between associativity and regularity of $\newprodinline$,
something we do not consider further here.

\section{The Principle of Local Reflexivity}

The classical Principle of Local Reflexivity states that for
a Banach space $E$, the local (or finite-dimensional)
structure of $E''$ is the same as that of $E$, taking account
of duality.  Formally, we have:

\begin{definition}
Let $E$ and $F$ be Banach spaces, and let $T\in\mc B(E,F)$.
For $\epsilon>0$, we say that $T$ is 
\emph{$(1+\epsilon)$-isomorphism onto its range} if
$(1-\epsilon)\|x\| \leq \|T(x)\| \leq (1+\epsilon)\|x\|$
for each $x\in E$.
\end{definition}

\begin{theorem}\label{PLR}
Let $E$ be a Banach space, and let $M\subseteq E''$ and $N\subseteq E'$
be finite-dimensional subspaces.  For each $\epsilon>0$ there exists
a $(1+\epsilon)$-isomorphism onto its range
$T:M\rightarrow E$ such that:
\begin{enumerate}
\item $\ip{\Phi}{\mu} = \ip{\mu}{T(\Phi)}$ for $\mu\in N$
   and $\Phi\in M$;
\item $T(\kappa_E(x)) = x$ for $x\in E$ such that
   $\kappa_E(x) \in M$.
\end{enumerate}
\end{theorem}
\begin{proof}
See \cite[Section~5.5]{Ryan} for a readable account.
\end{proof}

We wish to extend this result, using the results of Behrends in
\cite{Beh}.  Before we can do this, we need a word about tensor
products of Banach spaces.  Let $E$ and $F$ be Banach spaces,
and let $E\otimes F$ be the algebraic tensor product of $E$ with
$F$.  We define the \emph{projective tensor norm} on $E\otimes F$ by
\[ \|u\|_\pi = \inf\Big\{\sum_{i=1}^n \|x_i\| \|y_i\| : 
u = \sum_{i=1}^n x_i\otimes y_i \Big\}
\qquad (u\in E\otimes F). \]
Then the completion of $E\otimes F$ with respect to $\|\cdot\|_\pi$ is
the \emph{projective tensor product}, $E\proten F$.

See \cite{Ryan}, \cite{DF} or \cite[Section~VIII]{DU} for further details.
In particular, when $E$ is a
Banach space and $M$ is a \emph{finite-dimensional} Banach space,
then the dual of $\mc B(M,E)$ may be identified with $M\proten E'$ by
\[ \ip{x\otimes\mu}{T} = \ip{\mu}{T(x)}
\qquad (x\otimes\mu\in M\proten E', T\in\mc B(M,E)). \]
For general Banach spaces $F$ and $G$, the dual of $F\proten G$ is
$\mc B(F,G')$, under the identification
\[ \ip{T}{x\otimes y} = \ip{T(x)}{y}
\qquad (x\otimes y\in F\otimes G, T\in\mc B(F,G')). \]
Thus $\mc B(M,E)'' = (M\proten E')' = \mc B(M,E'')$, and we can check
that the canonical map $\kappa_{\mc B(M,E)}: \mc B(M,E)\rightarrow
\mc B(M,E'')$ satisfies $\kappa_{\mc B(M,E)}(T) = \kappa_E\circ T$
for $T\in\mc B(M,E)$.

The following definitions are from \cite{Beh}.  For a Banach space
$E$, we write $\FIN(E)$ for the collection of finite-dimensional
subspaces of $E$.

\begin{definition}
Let $E$ be a Banach space, and let $M\in\FIN(E'')$ and $N\in\FIN(E')$.
A map $T:M\rightarrow E$ is an $\epsilon$-\emph{isomorphism along $N$}
if $T$ is a $(1+\epsilon)$-isomorphism onto its range such that
$\ip{\Phi}{\mu} = \ip{\mu}{T(\Phi)}$ for $\mu\in N$ and $\Phi\in M$.

Let $(F_i)_{i=1}^n$ and $(G_j)_{j=1}^m$ be families of Banach
spaces.  Let $A_i:\mc B(M,E)\rightarrow F_i$ be an operator,
for $1\leq i\leq n$, and let $\psi_j:\mc B(M,E)\rightarrow G_j$
be an operator, for $1\leq j\leq m$.
For $1\leq i\leq n$, let $f_i\in F_i$, and for
$1\leq j\leq m$, let $C_j\subseteq G_j$ be a convex set.
Then $M$ satisfies:
\begin{enumerate}
\item the \emph{exact conditions} $(A_i,f_i)$, for
   $1\leq i\leq n$, and
\item the \emph{approximate conditions} $(\psi_j,C_j)$,
   for $1\leq j\leq m$,
\end{enumerate}
if for each $N\in\FIN(E')$ and $\epsilon>0$,
there exists an $\epsilon$-isomorphism $T:M\rightarrow E$ along $N$ such that
$A_i(T)=f_i$, for $1\leq i\leq n$, and $\psi_j(T)\in (C_j)_\epsilon
= \{ y+z : y\in C_j, z\in G_j, \|z\|\leq\epsilon \}$, for
$1\leq j\leq m$.
\end{definition}

Given $(A_i)$ as in the above definition, notice that we have that
$A_i' : F_i' \rightarrow M \proten E'$, and that $A_i'':
\mc B(M,E'') \rightarrow F_i''$.  For $M\subseteq E''$, let
$\iota_M:M\rightarrow E''$ be the inclusion map.

\begin{theorem}\label{behrends_result}
Let $E$ be a Banach space, $M\in\FIN(E'')$, and let $(F_i)$,
$(A_i)$, $(y_i)$, $(G_j)$, $(\psi_j)$ and
$(C_j)$ be as in the above definition.  Then the following are equivalent:
\begin{enumerate}
\item $M$ satisfies the exact conditions $(A_i,y_i)_{i=1}^n$
   and the approximate conditions\\  
   $(\psi_j,C_j)_{j=1}^m$;
\item $\iota_M$ is weak$^*$-continuous on the weak$^*$-closure
   of $A'_1(F'_1) + \cdots + A'_n(F'_n)$, $A''_i(\iota_M) =
   \kappa_{F_i}(y_i)$ for each $i$, and $\psi''_j(\iota_M)$ lies in
   the weak$^*$-closure of $\kappa_{G_j}(C_j)$, for each $j$.
\end{enumerate}

Suppose that the map $T\mapsto (A_i(T))_{i=1}^n$ from
$\mc B(M,E)$ to $A_1 \oplus \cdots\oplus A_n$ has a closed
range.  Then we may replace $\iota_M$ being weak$^*$-continuous
on the weak$^*$-closure of $\sum_{i=1}^n A'_i(F'_i)$ by there
existing some $T:M\rightarrow E$ which satisfies $A_i(T)=y_i$,
for $1\leq i\leq n$ ($T$ need not satisfy any other condition).
\end{theorem}
\begin{proof}
This is \cite[Theorem~2.3]{Beh}, and the remark thereafter.
\end{proof}

For example, let $A_M:\mc B(M,E) \rightarrow \mc B(M\cap\kappa_E(E),E)$
be the restriction operator, and let $B_M\in\mc B(M\cap\kappa_E(E),E)$
be the map $B_M(\kappa_E(x))=x$.  Then the Principle of
Local Reflexivity is just the statement that each $M\in\FIN(E'')$
satisfies the exact condition $(A_M,B_M)$.  Notice that
$A_M'':\mc B(M,E'')\rightarrow\mc B(M\cap\kappa_E(E),E'')$ is
also the restriction operator, and that $B_M'':M\cap\kappa_E(E)
\rightarrow E'$ is the inclusion map, so that condition (2) above
is easily verified in this case (or one can use the remark).

\section{Ultrapowers of modules}\label{mod_ultra}

We wish to extend the principle of local reflexivity to (bi)modules
of Banach algebras.  Let $\mc A$ be a Banach algebra, and let $E$ be
a Banach left $\mc A$-module (or a right $\mc A$-module, or an
$\mc A$-bimodule), so that
we can certainly apply the principle of local reflexivity to $E$.
However, we also want to take account of the $\mc A$-module structure,
that is, ensure that $T:M\rightarrow E$ is ``in some sense'' an
$\mc A$-module homomorphism (of course, $M$ will in general not
be a submodule).

It will be helpful to recall that $\kappa_E$ is an $\mc A$-module
homomorphism.  For the following, note that for $L\subseteq\mc A$ and
$M\subseteq E''$ finite-dimensional, we have that
\[ L\cdot M = \{ a\cdot\Phi : a\in L,\Phi\in M \} \in\FIN(E'').\]

\begin{theorem}\label{PLR_Mod_LeftRight}
Let $\mc A$ be a Banach algebra and let $E$ be a Banach left
$\mc A$-module.  Let $M\subseteq E''$, $L\subseteq\mc A$ and $N\subseteq E'$
be finite-dimensional, and let $\epsilon>0$.
Let $M_0\in\FIN(E'')$ be such that $L\cdot M + M\subseteq M_0$.
Then there exists $T:M_0\rightarrow E$, a
$(1+\epsilon)$-isomorphism onto its range, such that:
\begin{enumerate}
\item $\ip{\Phi}{\mu} = \ip{\mu}{T(\Phi)}$ for
   $\Phi\in M_0$ and $\mu\in N$;
\item $T(\kappa_E(x)) = x$ for $\kappa_E(x)
   \in M_0\cap \kappa_E(E)$;
\item $\| a\cdot T(\Phi) - T(a\cdot\Phi) \| \leq \epsilon\|a\|\|\Phi\|$
   for $a\in L$ and $\Phi\in M$.
\end{enumerate}
A similar result holds for Banach right $\mc A$-modules and
Banach $\mc A$-bimodules with condition (3) changed in the obvious way.
\end{theorem}
\begin{proof}
Let $\delta=\epsilon/5$ or $1$, whichever is smaller.
Let $(a_i)_{i=1}^n$ be a set in $L$ such that $\|a_i\|=1$ for
each $i$, and such that
\[ \min_{1\leq i\leq n} \|a-a_i\|<\delta
\qquad (a\in L, \|a\|=1). \]
For $1\leq i\leq n$, define $\psi_i : \mc B(M_0,E)\rightarrow
\mc B(M,E)$ by
\[ \psi_i(T)(\Phi) = T(a_i\cdot\Phi) - a_i\cdot T(\Phi)
\qquad (T\in\mc B(M_0,E), \Phi\in M). \]
Then $\psi'_i:M\proten E'\rightarrow M_0\proten E'$,
and for $\Phi\in M$, $\mu\in E'$ and $T\in\mc B(M_0,E)$,
we have
\[ \ip{\psi_i'(\Phi\otimes\mu)}{T} =
\ip{\mu}{\psi_i(T)(\Phi)} =
\ip{\mu}{T(a_i\cdot\Phi)-a_i\cdot T(\Phi)}, \]
so that $\psi_i'(\Phi\otimes\mu) = a_i\cdot\Phi\otimes\mu
- \Phi\otimes\mu\cdot a_i$.  Then, for $\Phi\in M$ and $\mu\in E'$,
we have
\[ \ip{\psi''_i(\iota_{M_0})}{\Phi\otimes\mu}
= \ip{\iota_{M_0}}{a_i\cdot\Phi\otimes\mu-\Phi\otimes\mu\cdot a_i}
= \ip{a_i\cdot\Phi}{\mu} - \ip{\Phi}{\mu\cdot a_i} = 0, \]
so that $\psi''_i(\iota_{M_0}) = 0$.

Consider the exact condition $(A_{M_0}, B_{M_0})$, as after
Theorem~\ref{behrends_result}.  Then clearly we have verified condition
(2) for the approximate conditions $(\psi_i,\{0\})$, for $1\leq i\leq n$.
Applying Theorem~\ref{behrends_result}, we find
$T\in\mc B(M_0,E)$, a $(1+\delta)$-isomorphism onto its
range, with conditions (1) and (2), and such that
$\|\psi_i(T)\|<\delta$ for $1\leq i\leq n$.
Then, for $a\in L$ and $\Phi\in M$ with
$\|a\|=\|\Phi\|=1$, we can find $i$ with $\|a-a_i\|<\delta$.
Then we have
\begin{align*}
\|a\cdot T(\Phi)-T(a\cdot\Phi)\| & \\
&\hspace{-10ex}\leq \|(a-a_i)\cdot T(\Phi)\|
+ \|a_i\cdot T(\Phi)-T(a_i\cdot\Phi)\| + \|T(a_i\cdot\Phi-a.\Phi)\| \\
&\hspace{-10ex}< \delta(1+\delta) + \|\psi_i(T)\| + (1+\delta)\delta
< 3\delta + 2\delta^2 < \epsilon.
\end{align*}
Thus we are done, as $\delta<\epsilon$.

Similarly, we can easily adapt the above argument to give the
result for right $\mc A$-modules and $\mc A$-bimodules.
\end{proof}

It would be nice if we could work with the exact conditions $(\phi_i,\{0\})$
above, but it is far from clear that we can apply condition (2)
of Theorem~\ref{behrends_result} in this case.  Fortunately, the
above is enough for our application.

We now apply this to prove a result about ultrapowers
of modules.  Notice that for a Banach algebra $\mc A$, a
left $\mc A$-module $E$ and an ultrafilter $\mc U$, we have
that $(E)_{\mc U}$ is a left $\mc A$-module with pointwise
module action.  As usual, $E'$ becomes a right $\mc A$-module
for the module action given by
\[ \ip{\mu\cdot a}{x} = \ip{\mu}{a\cdot x}
\qquad (a\in\mc A, x\in E, \mu\in E'). \]
Similarly, $E''$ becomes a left $\mc A$-module.
Recall the map $\sigma_{\mc U}:(E)_{\mc U}
\rightarrow E''$; this is easily seen to be a left
$\mc A$-module homomorphism, as
\[ \ip{a \cdot \sigma_{\mc U}(x)}{\mu} =
\lim_{i\rightarrow\mc U} \ip{\mu\cdot a}{x_i}
= \lim_{i\rightarrow\mc U} \ip{\mu}{a\cdot x_i}
= \ip{\sigma_{\mc U}(a\cdot x)}{\mu}, \]
for $a\in\mc A,\mu\in E'$ and $x=(x_i)\in (E)_{\mc U}$.
This all holds, with obvious modifications, for right
$\mc A$-modules and $\mc A$-bimodules.

\begin{theorem}\label{mod_ultra_thm}
Let $\mc A$ be a Banach algebra, and let $E$ be a
left $\mc A$-module.  Then there exists an ultrafilter
$\mc U$ and a map $K:E''\rightarrow (E)_{\mc U}$ such that:
\begin{enumerate}
\item $K$ is an isometry and a left $\mc A$-module homomorphism;
\item $\sigma_{\mc U} \circ K$ is the identity on $E''$;
\item $K\circ\kappa_E:E\rightarrow (E)_{\mc U}$ is the
   canonical map.
\end{enumerate}
Similar results hold for right $\mc A$-modules and
$\mc A$-bimodules.
\end{theorem}
\begin{proof}
We carefully prove this sort of result once, for completeness.
Define
\[ I = \big\{ (M,N,L,\epsilon) : M\in\FIN(E''), N\in\FIN(E'),
L\in\FIN(\mc A), \epsilon>0 \big \} \]
with partial order given by $(M_1,N_1,L_1,\epsilon_1) \geq
(M_2,N_2,L_2,\epsilon_2)$ if and only if $M_1 \supseteq M_2$,
$N_1 \supseteq N_2$, $L_1\supseteq L_2$, and
$\epsilon_1 \leq \epsilon_2$.  Then $(I,\leq)$ becomes
a directed set, and we let $\mc U$ be some ultrafilter
refining the order-filter on $I$.

We now define $K:E''\rightarrow(E)_{\mc U}$.  Fix $\Phi\in E''$,
and let $K(\Phi)=(x_i)_{i\in I} \in (E)_{\mc U}$, where, for
$i=(M,N,L,\epsilon)\in I$, we let $x_i = T(\Phi)$ where $T:
M\rightarrow E$ is given by the above theorem (with, say,
$M_0 = L\cdot M+M$).  Then, if $\Phi,\Psi\in E''$, let
$K(\Phi) = (x_i)$, $K(\Psi) = (y_i)$ and $K(\Phi+\Psi)=(z_i)$,
so that for $i=(M,N,L,\epsilon)$ with $\Phi,\Psi\in M$,
we have that $x_i + y_i = T(\Phi)+T(\Psi) = T(\Phi+\Psi)
= z_i$.  Thus $(x_i) + (y_i) = (z_i)$ in $(E)_{\mc U}$, so that
$K$ is linear.  Similarly, we see that $K$ is an isometry,
as $\epsilon\rightarrow 0$ along $\mc U$.  It follow from conditions
(1) and (2) in the above theorem, that, respectively,
$\sigma_{\mc U} \circ K$ is the identity on $E''$, and
$K\circ\kappa_E$ is the canonical map $E\rightarrow (E)_{\mc U}$.

Finally, we show that $K$ is a left $\mc A$-module homomorphism.
Let $a\in\mc A$ and $\Phi\in E''$, let $K(\Phi)=(x_i)$ and
$K(a\cdot\Phi)=(y_i)$, and let $i=(M,N,L,\epsilon)$ with
$\Phi,a\cdot\Phi\in M$ and $a\in L$.  Then
\[ \| a\cdot x_i - y_i \| = \| a\cdot T(\Phi) - T(a\cdot\Phi) \|
\leq \epsilon\|a\|\|\Phi\|, \]
by condition (3) in the above theorem.  Thus $a\cdot K(\Phi) =
K(a\cdot\Phi)$ in $(E)_{\mc U}$.
\end{proof}

\section{Ultrapowers of algebras}

Let $\mc A$ be a Banach algebra.  Using the above, as $\mc A''$
is an $\mc A$-bimodule, we can find an ultrafilter $\mc U$ and
an isometry $K:\mc A''\rightarrow(\mc A)_{\mc U}$ which is an
$\mc A$-module homomorphism, and with $\sigma_{\mc U}\circ K$
the identity on $\mc A''$.  This is not quite enough to show
that $\newprodinline$ agrees with either Arens product on $\mc A''$.
We now show that, at least when $\mc A$ is Arens regular,
we can do better.

\begin{theorem}\label{PLR_alg}
Let $\mc A$ be an Arens regular Banach algebra,
$M\in\FIN(\mc A'')$, $N\in\FIN(\mc A')$ and $\epsilon>0$.
Let $M_0=M+M\aone M$ and $N_0=N+M\cdot N$.  Then there exists
a $(1+\epsilon)$-isomorphism onto its range
$T:M_0\rightarrow\mc A$ such that:
\begin{enumerate}
\item $\ip{\Phi}{\mu}=\ip{\mu}{T(\Phi)}$ for
   $\Phi\in M_0$ and $\mu\in N_0$;
\item $T(\kappa_{\mc A}(a))=a$ for $\kappa_{\mc A}(a)
   \in M_0\cap\kappa_{\mc A}(\mc A)$;
\item $|\ip{\mu}{T(\Phi\aone\Psi)-T(\Phi)T(\Psi)}|
   \leq \epsilon\|\mu\|\|\Phi\|\|\Psi\|$ for
   $\mu\in N$ and $\Phi,\Psi\in M$.
\end{enumerate}
\end{theorem}
\begin{proof}
Let $\delta>0$ be such that $\delta<\epsilon$ and
$\delta(1+\delta)(3+\delta)<\epsilon$.
Let $(\mu_i)_{i=1}^n\subseteq N$ be
such that $\|\mu_i\|=1$ for each $i$, and such that we have
\[ \min_{1\leq i\leq n} \|\mu_i-\mu\| < \delta
\qquad (\mu\in N, \|\mu\|=1). \]
For $1\leq i\leq n$, define $\psi_i : \mc B(M_0,\mc A)
\rightarrow \mc B(M_0,\mc A')$ by
\[ \psi_i(T)(\Phi) = T(\Phi)\cdot\mu_i
\qquad (T\in\mc B(M_0,\mc A), \Phi\in M_0), \]
and define $T_i\in\mc B(M_0,\mc A')$ by
$T_i(\Phi)=\Phi\cdot\mu_i$ for $\Phi\in M_0$.
Then we have $\psi'_i:M_0\proten\mc A''\rightarrow
M_0\proten\mc A'$, and, for $T\in\mc B(M_0,\mc A)$, $\Phi\in M_0$
and $\Lambda\in\mc A''$, we have
\[ \ip{\psi'_i(\Phi\otimes\Lambda)}{T} =
\ip{\Lambda}{\psi_i(T)(\Phi)} =
\ip{\Lambda}{T(\Phi)\cdot\mu_i} =
\ip{\mu_i\cdot\Lambda}{T(\Phi)}. \]
Thus we have $\psi'_i(\Phi\otimes\Lambda) =
\Phi \otimes \mu_i\cdot\Lambda$, and so
\[ \ip{\psi''_i(\iota_{M_0})}{\Phi\otimes\Lambda}
= \ip{\Phi}{\mu_i\cdot\Lambda} = \ip{\Lambda\atwo\Phi}{\mu_i}
= \ip{\Lambda\aone\Phi}{\mu_i} = \ip{\Lambda}{\Phi\cdot\mu_i} =
\ip{\kappa_{\mc A'}(T_i(\Phi))}{\Lambda}, \]
as $\mc A$ is Arens regular.  Thus $\psi''_i(\iota_{M_0})
= \kappa_{\mc B(M_0,\mc A')}(T_i)$.
Again, we can then find $T\in\mc B(M_0,\mc A)$ satisfying
(1) and (2), and such that $\|\psi_i(T)-T_i\|<\delta$ for
$1\leq i\leq n$.

For $\mu\in N$ and $\Phi,\Psi\in M$ with
$\|\mu\|=\|\Phi\|=\|\Psi\|=1$, let $i$ be such that $\|\mu-\mu_i\|
<\delta$.  Then $\Phi\aone\Psi\in M_0$ and $\Psi\cdot\mu\in N_0$
so that we have
\[ \ip{\mu}{T(\Phi\aone\Psi)} = \ip{\Phi\aone\Psi}{\mu}
= \ip{\Phi}{\Psi\cdot\mu} = \ip{\Psi\cdot\mu}{T(\Phi)}. \]
As $\|\psi_i(T)-T_i\|<\delta$, we have
$\| T(\Psi)\cdot\mu_i - \Psi\cdot\mu_i \| < \delta$, and so
\begin{align*}
\|T(\Psi)\cdot\mu-\Psi\cdot\mu\| & \\
&\hspace{-10ex}\leq
\|T(\Psi)\cdot\mu - T(\Psi)\cdot\mu_i\| + \|T(\Psi)\cdot\mu_i-
\Psi\cdot\mu_i\| + \|\Psi\cdot\mu_i-\Psi\cdot\mu\| \\
&\hspace{-10ex}< \delta\|T(\Psi)\| + \delta + \delta\|\Psi\|
\leq \delta(1+\delta) + 2\delta.
\end{align*}
Putting these together, we then get
\begin{align*}
|\ip{\mu}{T(\Phi\aone\Psi)-T(\Phi)T(\Psi)}| &=
|\ip{\Psi\cdot\mu}{T(\Phi)} - \ip{T(\Psi)\cdot\mu}{T(\Phi)}| \\
&< \|T(\Phi)\| \big( \delta(1+\delta) + 2\delta \big)
\leq (1+\delta)\big( \delta(1+\delta) + 2\delta \big) \\
&= \delta(1+\delta)(3+\delta) < \epsilon,
\end{align*}
as required, completing the proof.
\end{proof}

\begin{theorem}\label{arens_case}
Let $\mc A$ be an Arens regular Banach algebra.  There exists an
ultrafilter $\mc U$ and an isometry $K:\mc A''\rightarrow
(\mc A)_{\mc U}$ such that:
\begin{enumerate}
\item $\sigma_{\mc U}\circ K$ is the identity on $A''$;
\item $K\circ\kappa_{\mc A}$ is the canonical map $\mc A
   \rightarrow (\mc A)_{\mc U}$;
\item $\newprodinline$, defined using $K$, agrees with the Arens
   products on $\mc A''$.
\end{enumerate}
\end{theorem}
\begin{proof}
This follows exactly as in the proof of
Theorem~\ref{mod_ultra_thm}.
\end{proof}

Let $\mc A$ be an Arens regular Banach algebra, and form
$K:\mc A''\rightarrow(\mc A)_{\mc U}$ as in the theorem.  For
$\Phi,\Psi\in\mc A''$, let $(a_i)=K(\Phi)$ and $(b_i)=K(\Psi)$,
so that for $\mu\in\mc A'$,
\begin{gather*}
\ip{\Phi}{\mu} = \lim_{i\rightarrow\mc U} \ip{\mu}{a_i}, \quad
\ip{\Psi}{\mu} = \lim_{i\rightarrow\mc U} \ip{\mu}{b_i}, \\
\ip{\Phi\aone\Psi}{\mu} = \ip{\Phi\newprod\Psi}{\mu}
   = \lim_{i\rightarrow\mc U} \ip{\mu}{a_ib_i}.
\end{gather*}
Compare this symmetric definition of $\Phi\aone\Psi$ to
the formulae in Section~\ref{Intro}.  Of course, here we have to
be careful in our choice of $K$.

As in the introduction, we note that \cite[Theorem~12]{IL} shows that
it is too much to expect the above to be true for a non-Arens regular
Banach algebra $\mc A$, with (3) replaced by asking for $\newprodinline$
to agree with $\aone$ or $\atwo$.  At least, this is true if $\mc A$
is commutative.  It would be interesting to know if we could ever
have, say, $\newprodinline = \aone$ for a non-commutative, non-Arens
regular Banach algebra.

\subsection{Asking for an algebra homomorphism}

A much stronger result than the above would be to find $\mc U$ and
$K:\mc A''\rightarrow(\mc A)_{\mc U}$ with $K$ being an algebra
homomorphism (presumably assuming that $\mc A$ is Arens
regular).  We now present a case when this is possible.

\begin{proposition}
Let $\mc A$ be a commutative, Arens regular Banach algebra such that
$\mc A$ is an essential ideal in $\mc A''$.  Suppose that $\mc A$ has an
approximate identity $(e_\alpha)$ consisting of idempotents, which is
bounded in the multiplier norm.
Then there exists an ultrafilter $\mc U$ and $K:\mc A''\rightarrow
(\mc A)_{\mc U}$ with $K$ being an algebra homomorphism, and such
that $\sigma\circ K$ is the identity on $\mc A''$, and $K\circ\kappa_{\mc A}$
if the canonical inclusion $\mc A\rightarrow(\mc A)_{\mc U}$.
\end{proposition}
\begin{proof}
Let $(e_\alpha)$ be indexed by the directed set $I$, and let $\mc U$
be an ultrafilter on $I$ refining the order filter.  As $\mc A$ is an
ideal in $\mc A''$, we see that for $\Phi\in\mc A''$, we have that
$\Phi e_\alpha \in \mc A$ for each $\alpha\in I$.  As $(e_\alpha)$ is
bounded in the multiplier norm, there exists $M>0$ such that
\[ \|a e_\alpha\| \leq M \|a\| \qquad (a\in\mc A, \alpha\in I). \]
Thus $\|e_\alpha\cdot\mu\|\leq M\|\mu\|$ for $\mu\in\mc A'$ and each
$\alpha$, and so $\|\Phi e_\alpha\|\leq M\|\Phi\|$ for $\Phi\in\mc A''$
and each $\alpha$.  Hence we may define $K$ by
\[ K(\Phi) = (\Phi e_\alpha) \in (\mc A)_{\mc U} \qquad (\Phi\in\mc A''). \]
$K$ is hence linear and bounded, and as $(e_\alpha)$ is an approximate
identity, we see that $K\circ\kappa_{\mc A}$ is the canonical map
$\mc A\rightarrow(\mc A)_{\mc U}$.  For each $\alpha$, we know that
$e_\alpha^2=e_\alpha$, and so
\[ K(\Phi)K(\Psi) = (\Phi e_\alpha \Psi e_\alpha)
= (\Phi\Psi e_\alpha e_\alpha) = K(\Phi\Psi) \qquad (\Phi,\Psi\in\mc A''). \]
For $\mu\in\mc A'$ and $a\in\mc A$, we see that
\[ \ip{\sigma K(\Phi)}{a\cdot\mu}
= \lim_{\alpha\rightarrow\mc U} \ip{\Phi}{e_\alpha a\cdot\mu}
= \ip{\Phi}{a\cdot\mu}. \qquad (\Phi\in\mc A''). \]
Hence $\sigma K(\Phi)=\Phi$ when restricted to $\mc A\cdot\mc A'$.
If $\Phi \in (\mc A\cdot\mc A')^\perp$, then $\Phi a=0$ for each
$a\in\mc A$.  As $\mc A$ is an essential ideal in $\mc A''$, by
definition, $\Phi=0$, and so we conclude that $\sigma\circ K$
is the identity on $\mc A''$.
\end{proof}

The above applies in particular to $\mc A=c_0$.  It would be
interesting to know if the same conclusions hold in the
\emph{non-commutative} case, namely the compact operators on $\ell^2$.

Following an example from \cite{GH}, let
$\mc A = C([0,1])$, a commutative C$^*$-algebra.  Then $\mc A''$ is a von
Neumann algebra, and so has many (self-adjoint) projections.
Let $\mc U$ be an ultrafilter, and suppose that $\mc A''$ is isomorphic, as
a Banach algebra, to a subalgebra of $(\mc A)_{\mc U}$.  Thus $(\mc A)_{\mc U}$
contains non-trivial, not necessarily self-adjoint, projections.  However,
$(\mc A)_{\mc U}$ is isomorphic to $C(K)$ for some compact Hausdorff $K$,
and any projection in $C(K)$ is automatically self-adjoint.  It follows from
\cite[Proposition~2.1]{GH} that $(\mc A)_{\mc U}$ does not contain non-trivial
projections, as the only projections in $\mc A$ are $0$ and $1$, and so we
have a contradiction.  It would be interesting to, say, characterise which
C$^*$-algebras $\mc A$ are such that $\mc A''$ is isomorphic, or $*$-isomorphic,
to a subalgebra of an ultrapower of $\mc A$.  This of course has links
to the notorious Connes-embedding problem for von Neumann algebras.
 
Following \cite[Example~III.1]{GI}, let $\mc A=\ell^1$ with the
pointwise product.  Then $c_0 \subseteq (\ell^1)'=\ell^\infty$, and we
can decompose $(\ell^1)''$ as $c_0^\perp \oplus \ell^1$, where
\[ c_0^\perp = \{ \Phi\in(\ell^1)'' : \ip{\Phi}{x}=0 \ (x\in c_0) \}. \]
Furthermore, $\|\Phi+a\| = \|\Phi\|+\|a\|$ for $\Phi\in c_0^\perp$ and
$a\in\ell^1$.  Then the product on $(\ell^1)''$ is simply
\[ (\Phi,a) (\Psi,b) = ab \qquad (a,b\in\ell^1, \Phi,\Psi\in c_0^\perp). \]
Thus $\mc A$ is an ideal in $\mc A''$, and clearly $\mc A$ has an approximate
identity consisting of idempotents, and bounded in the multiplier norm.
Notice that $\mc A$ is certainly not an essential ideal in $\mc A''$.

\begin{proposition}
With notation as above, there exists a ultrafilter $\mc U$ and an
isometry $K:(\ell^1)'' \rightarrow (\ell^1)_{\mc U}$ satisfying the
conclusions of Theorem~\ref{arens_case}, and with $K$ being a homomorphism.
\end{proposition}
\begin{proof}
Let $M\in\FIN((\ell^1)'')$, $N\in\FIN(\ell^\infty)$ and $\epsilon>0$.
As $(\ell^1)''$ is an $\mathcal{L}_1$ space (see \cite[Chapter~2]{Ryan},
for example), by enlarging $M$ is
necessary, we can find a basis $(m_i)_{i=1}^k$ for $M$ such that
$\|m_i\|=1$ for each $i$, and
\[ (1-\epsilon)\sum_{i=1}^k |\alpha_i| \leq
\Big\| \sum_{i=1}^k \alpha_i m_i \Big\|
\leq (1+\epsilon)\sum_{i=1}^k |\alpha_i|
\qquad ( (\alpha_i)_{i=1}^k\subseteq\mathbb C ). \]
As $(\ell^1)'' = c_0^\perp \oplus \ell^1$, we may suppose that
$m_i\in\ell^1$ for $1\leq i< \hat k$ and $m_i\in c_0^\perp$ for
$\hat k\leq i\leq k$.
Similarly, it is no loss of generality to suppose that $N$ is
the linear span of indicator functions $\chi_{A_1},\ldots,\chi_{A_l}$
with $(A_i)_{i=1}^l \subseteq\mathbb N$ pairwise disjoint.  Furthermore,
we may suppose that
\[ \sup\big\{ |\ip{m}{x}| : x\in N, \|x\|\leq 1 \big\}
\geq (1-\epsilon)\|m\| \qquad (m\in M). \]
Order $(A_i)$ so that $A_i$ is finite for $1\leq i<\hat l$, and
$A_i$ is infinite otherwise.

Let $a_i=m_i\in\ell^1$ for $1\leq i<\hat k$, let
$a_i = (a^{(i)}_t)_{t\in\mathbb N}\in\ell^1$.  Suppose we have chosen
$a_1,\ldots,a_r$, and let $N$ be such that $\sum_{|t|>N} |a^{(i)}_t|
< \epsilon$ for $1\leq i\leq r$, and such that $A_j\subseteq
\{1,\cdots,N\}$ for $1\leq j<\hat l$.  As $m_{r+1}\in c_0^\perp$, for
all finite sets $A\subseteq\mathbb N$, we have that
\[ \ip{m_{r+1}}{\chi_{A_j}} = 0, \ (1\leq j<\hat l), \quad
\ip{m_{r+1}}{\chi_{A_j\setminus A}} = \ip{m_{r+1}}{\chi_{A_j}}, \
(\hat l\leq j\leq l). \]
Consequently, by a simple argument, we can find $a_{r+1}=(a^{(r+1)}_t)\in\ell^1$
with $\|a_{r+1}\|=\|m_{r+1}\|=1$, $\ip{\chi_{A_j}}{a_{r+1}} =
\ip{m_{r+1}}{\chi_{A_j}}$ for $\hat l\leq j\leq l$, with $a^{(r+1)}_t=0$ for $t\leq N$,
and with $|a^{(r+1)}_t|<\epsilon$ for all $t$.
Thus $\ip{\chi_{A_j}}{a_{r+1}} = 0 = \ip{m_{r+1}}{\chi_{A_j}}$ for
$1\leq j<\hat l$.

We hence find $(a_i)_{i=1}^k\in\ell^1$ such that $\ip{x}{a_i} = \ip{m_i}{x}$
for $x\in N$ and $1\leq i\leq k$.  Furthermore, $(a_i)_{i=\hat k}^k$ have
disjoint support.  Consequently, for $(\alpha_i)_{i=1}^k\subseteq\mathbb C$,
\[ \Big\| \sum_{i=1}^k \alpha_i a_i \Big\|
= \Big\| \sum_{i=1}^{\hat k-1} \alpha_i m_i \Big\|
   + \sum_{i=\hat k}^k |\alpha_i|
\leq (1+\epsilon) \sum_{i=1}^k |\alpha_i|
\leq \frac{1+\epsilon}{1-\epsilon} \Big\| \sum_{i=1}^k \alpha_i m_i \Big\|. \]
Similarly,
\[ \Big\| \sum_{i=1}^k \alpha_i a_i \Big\|
\geq \sup\Big\{ \sum_{i=1}^k \alpha_i \ip{x}{a_i} : x\in N, \|x\|=1 \Big\}
\geq (1-\epsilon) \Big\| \sum_{i=1}^k \alpha_i m_i \Big\|. \]
Hence the map $T:M\rightarrow\ell^1$ defined by $T(m_i)=a_i$ and linearity
is a $(1+\hat\epsilon)$-isomorphism onto its range, for $\hat\epsilon = 
2\epsilon(1-\epsilon)^{-1}$.  Furthermore, $\ip{m}{x} = \ip{x}{T(m)}$
for $m\in M$ and $x\in N$, and $T\kappa_{\ell^1}(a)=a$ for $a\in\ell^1$
with $\kappa_{\ell^1}(a)\in M$.

By a now standard argument, we find an isometry $K:(\ell^1)'' \rightarrow
(\ell^1)_{\mc U}$ such that $\sigma_{\mc U}K$ is the identity on $(\ell^1)''$,
$K\kappa_{\ell^1}$ is the canonical map $\ell^1\rightarrow (\ell^1)_{\mc U}$,
and for $\Phi,\Psi\in c_0^\perp$, we have that $K(\Phi)K(\Psi)=0$.
This last fact follows as we chose $(a_i)_{i=\hat k}^k$ above with disjoint
support, and with $\|a_i\|_\infty$ small for $i\geq\hat k$ (this deals with
the case that $\Psi$ is a scalar multiple of $\Phi$).  By the discussion
above $K$ is hence an algebra homomorphism, as required.
\end{proof}

\section{Dual Banach algebras and weakly almost periodic functionals}
\label{wap_ap_sec}

Let $\mc A$ be a Banach algebra such that $\mc A = E'$ for some
Banach space $E$.  We say that $\mc A$ is a \emph{dual Banach algebra}
if the product on $\mc A$ is separately weak$^*$-continuous.  It is
shown in \cite[Section~1]{Runde} (see also \cite[Section~2]{DawsDBA})
that $\mc A$ is a dual Banach algebra if and only if
$\kappa_E(E)\subseteq E''=\mc A'$ is an $\mc A$-submodule.

Notice that $\mc A''$ is always a dual Banach space.  It is not hard
to show that $\mc A''$, with either Arens product, is a dual Banach
algebra if and only if $\mc A$ is Arens regular.  Similarly, by
the definition in \cite{IL}, the product $\newprodinline$ is
regular if and only if it makes $\mc A''$ a dual Banach algebra
(of course, taking account of the fact that $\newprodinline$
may not be associative).  Hence \cite[Proposition~6]{IL} shows that
$(\mc A'',\newprodinline)$ is a dual Banach algebra (in this
not necessarily associative sense) if and only if $\newprodinline=
\aone=\atwo$ (in which case, $\newprodinline$ is automatically
associative).  In particular, $\newprodinline$ cannot turn $\mc A''$ into
a genuine dual Banach algebra unless $\mc A$ is already Arens regular.

Let $\mc A$ be a Banach algebra.  We say that $\mu\in\mc A'$ is
\emph{weakly almost periodic} if the map $\mc A\rightarrow\mc A';
a\mapsto a\cdot\mu$ is weakly compact.  We write $\mu\in\wap(\mc A')$
in this case (some authors write $\wap(\mc A)$ for this).  The
space $\wap(\mc A')$ has been widely studied, especially in the
context of group algebras.  We now collect some useful results.

\begin{proposition}
Let $\mc A$ be a Banach algebra.  Then $\wap(\mc A')$ is a
closed submodule of $\mc A'$.  For $\mu\in\mc A'$, we have that
$\mu\in\wap(\mc A')$ if and only if $\ip{\Phi\aone\Psi}{\mu} =
\ip{\Phi\atwo\Psi}{\mu}$ for $\Phi,\Psi\in\mc A''$.  In particular,
$\mc A$ is Arens regular if and only if $\wap(\mc A')=\mc A'$.

Let $X\subseteq\mc A'$ be a closed submodule, so we identify $X'$
with the quotient $\mc A''/X^\perp$.  The following are equivalent:
\begin{enumerate}
\item $X\subseteq\wap(\mc A')$;
\item the Arens products drop to a well-defined product on $X'$
turning $X'$ into a dual Banach algebra.
\end{enumerate}
\end{proposition}
\begin{proof}
These facts are collected in \cite[Section~2]{DawsDBA}.
The first result is due to John Pym.  The second result can be
found in many places in the literature: see \cite[Lemma~1.4]{LL2}
for example, which shows this for commutative Banach algebras.
\end{proof}

Similarly, we say that $\mu\in\mc A'$ is \emph{almost periodic},
written $\mu\in\ap(\mc A')$, if the map $\mc A\rightarrow\mc A';
a\mapsto a\cdot\mu$ is (norm) compact.  

\begin{proposition}\label{ap_char}
Let $\mc A$ be a Banach algebra.  Then $\ap(\mc A')$ is a
closed submodule of $\mc A'$.

Let $X\subseteq\mc A'$ be a closed submodule, so we identify $X'$
with the quotient $\mc A''/X^\perp$.  The following are equivalent:
\begin{enumerate}
\item $X\subseteq\ap(\mc A')$;
\item the Arens products drop to a well-defined product on $X'$
which is jointly continuous on bounded spheres.
\end{enumerate}
\end{proposition}
\begin{proof}
Lau shows this for a certain class of commutative Banach
algebras in \cite[Theorem~5.8]{Lau}, although
the proof is very easy to adapt to the general case.
Compare also \cite[Proposition~7]{IL}.
\end{proof}

Let $\mc A$ be a Banach algebra and let $\mc U$ be an ultrafilter.
Define $\sigma_{\mc U}^{{\wap}} : (\mc A)_{\mc U} \rightarrow \wap(\mc A')'$
by
\[ \ip{\sigma_{\mc U}^{{\wap}}((a_i))}{\mu} = \lim_{i\rightarrow\mc U}
\ip{\mu}{a_i} \qquad ((a_i)\in(\mc A)_{\mc U}, \mu\in\wap(\mc A')). \]
That is, $\sigma_{\mc U}^{{\wap}}$ is simply the map $\sigma_{\mc U}$
composed with the quotient map $\mc A''\rightarrow\wap(\mc A')'$.

\begin{theorem}\label{PLR_alg_wap}
Let $\mc A$ be a Banach algebra, let $M\in\FIN(\mc A'')$, $N\in\FIN(\wap(\mc A'))$
and $\hat N\in\FIN(\mc A')$, and let $\epsilon>0$.
Let $M_0=M+M\aone M$ and $N_0=N+M\cdot N$.  Then there exists
a $(1+\epsilon)$-isomorphism onto its range
$T:M_0\rightarrow\mc A$ such that:
\begin{enumerate}
\item $\ip{\Phi}{\mu}=\ip{\mu}{T(\Phi)}$ for
   $\Phi\in M_0$ and $\mu\in \hat N$;
\item $T(\kappa_{\mc A}(a))=a$ for $\kappa_{\mc A}(a)
   \in M_0\cap\kappa_{\mc A}(\mc A)$;
\item $|\ip{\mu}{T(\Phi\aone\Psi)-T(\Phi)T(\Psi)}|
   \leq \epsilon\|\mu\|\|\Phi\|\|\Psi\|$ for
   $\mu\in N$ and $\Phi,\Psi\in M$.
\end{enumerate}
\end{theorem}
\begin{proof}
If we examine the proof of Theorem~\ref{PLR_alg}, we see that it
will hold for non-Arens regular Banach algebras, so long as $N\subseteq
\wap(\mc A')$.  It is an easy exercise to take account of $\hat N$.
\end{proof}

\begin{corollary}
Let $\mc A$ be a Banach algebra.  There exists an
ultrafilter $\mc U$ and an isometry $K:\mc A''\rightarrow
(\mc A)_{\mc U}$ such that:
\begin{enumerate}
\item $\sigma_{\mc U}\circ K$ is the identity on $A''$;
\item $K\circ\kappa_{\mc A}$ is the canonical map $\mc A
   \rightarrow (\mc A)_{\mc U}$;
\item Let $\iota = \iota_{{\wap}}:\mc A''\rightarrow\wap(\mc A')'$
   be the quotient map.  For $\Phi,\Psi\in\mc A''$, we have that
   \[ \ip{\iota(\Phi) \iota(\Psi)}{\mu} =
   \ip{\sigma_{\mc U}^{\wap}( K(\Phi) K(\Psi) )}{\mu}
   \qquad (\mu\in\wap(\mc A')). \]
\end{enumerate}
\end{corollary}

It would seem to be natural to ask if we could define $K$ as a map $\wap(\mc A')'
\rightarrow (\mc A)_{\mc U}$ with $\sigma_{\mc U}^{{\wap}}\circ K$
the identity on $\wap(\mc A')'$.

\begin{proposition}
Let $\mc A$ be a Banach algebra.  The following are equivalent:
\begin{enumerate}
\item\label{wap_K_map} There exists a map
$K_{\wap}:\wap(\mc A')' \rightarrow (\mc A)_{\mc U}$
with $\sigma_{\mc U}^{\wap}\circ K_{\wap}$ the identity on $\wap(\mc A')'$;
\item\label{wap_wcom} $\wap(\mc A')^\perp$ is complemented in $\mc A''$ (that is,
$\wap(\mc A')$ is \emph{weakly-complemented});
\end{enumerate}
\end{proposition}
\begin{proof}
If (\ref{wap_K_map}) holds, then $L=\sigma_{\mc U}\circ K$ is a map
$\wap(\mc A')' \rightarrow \mc A''$.  Then $\iota_{\wap} \circ L$ is the
identity on $\wap(\mc A')'$, and so $L\circ\iota_{{\wap}}$ is a projection of
$\mc A''$ onto the image of $L$, with complementary space $\wap(\mc A')^\perp$,
and so (\ref{wap_wcom}) holds.

If (\ref{wap_wcom}) holds then let $P:\mc A''\rightarrow\wap(\mc A')^\perp$
be a projection.  Let $K:\mc A''\rightarrow(\mc A)_{\mc U}$ be such that
$\sigma_{\mc U}\circ K$ is the identity on $\mc A''$.  We identify $\wap(\mc A')'$
with the quotient $\mc A''/\wap(\mc A')^\perp$, and let $\iota_{\wap}$ be
the quotient map.  Define $K_{\wap}:\wap(\mc A')'\rightarrow(\mc A)_{\mc U}$ by
\[ K_{\wap}(\iota_{\wap}(\Phi)) = K\big(\Phi - P(\Phi)\big) \qquad (\Phi\in\mc A''). \]
Then, if $\iota_{\wap}(\Phi) = \iota_{\wap}(\Psi)$, we have that $\Phi-\Psi
\in\wap(\mc A')^\perp$, so that $P(\Phi-\Psi)=\Phi-\Psi$, and so
$K(\Phi-\Psi-P(\Phi-\Psi))=0$, showing that $K_{\wap}$ is well-defined.
Then $\sigma_{\mc U}^{\wap} K_{\wap} \iota_{\wap} = \iota_{\wap} \sigma_{\mc U} K (I-P)
= \iota_{\wap} (I-P) = \iota_{\wap}$, and so (\ref{wap_K_map}) holds.
\end{proof}

We currently have no examples showing that $\wap(\mc A')^\perp$ can be
complemented in $\mc A''$.  An obvious place to look is at group algebras
$L^1(G)$, for which $\wap(L^\infty(G))$ is (reasonably) well understood
(see \cite[Section~7]{DawsDBA} and references therein).

We can similarly define $\sigma_{\mc U}^{{\ap}}$, and as $\ap(\mc A')
\subseteq \wap(\mc A')$, all of the above holds for $\ap$, with suitable
modifications.  By Proposition~\ref{ap_char}, we see that
$\sigma_{\mc U}^{{\ap}}$ is actually an \emph{algebra homomorphism}.

\section{Automatic regularity for C$^*$-algebras}

As stated in the introduction, \cite[Corollary~II.2]{GI} shows that
when $\mc A$ is a C$^*$-algebra the map $\newprodinline$ always agrees with
$\aone=\atwo$, so long as $K:\mc A''\rightarrow(\mc A)_{\mc U}$ is an
\emph{isometry} (and satisfies our two standing assumptions, namely
that $\sigma\circ K$ is the identity on $\mc A''$, and that
$K\circ\kappa_{\mc A}$ is the canonical embedding).

The proof of \cite[Proposition~7]{IL} is trivially adapted to our
situation (that is, when $K$ is not assumed to be an isometry) to show
that when $\ap(\mc A')=\mc A'$, we have that $\newprodinline=\aone$
(as $\wap(\mc A')=\mc A'$, such algebras are always Arens regular).

Recall that a positive functional of a C$^*$-algebra is \emph{pure}
if the associated GNS representation is irreducible.
Recall that a C$^*$-algebra is said to be \emph{scattered} if
every positive functional is a sum of pure positive functionals.
For a commutative C$^*$-algebra, $C(X)$ is scattered if and only if
$X$ is scattered, in the topological sense that every subset of $X$
contains an isolated point.  Then \cite[Theorem~3.2]{Quigg} shows the following:

\begin{theorem}
For a C$^*$-algebra $\mc A$, $\ap(\mc A')=\mc A'$ if and only if
A is scattered and each irreducible representation
of $\mc A$ is finite-dimensional.
\end{theorem}

It seems possible that for well-behaved Banach algebras (say, certainly
for C$^*$-algebras), that when $\mc A' \not= \ap(\mc A')$, we can
construct a map $K:\mc A''\rightarrow(\mc A)_{\mc U}$ satisfying our
two assumptions, but with $\newprodinline \not= \aone$.  We can currently
only show this under an extra assumption.

Before we proceed, we provide a general way of finding such maps $K$.

\begin{proposition}\label{when_K_not_good}
Let $\mc A$ be an Arens regular Banach algebra.  Suppose that there
exist weakly-null nets $(a_\alpha)$ and $(b_\alpha)$, on the same index set,
such that $(a_\alpha b_\alpha)$ is not weakly-null.  Then there exists
an ultrafilter $\mc U$ and a bounded map $K:\mc A''\rightarrow(\mc A)_{\mc U}$
which satisfies our standing assumptions, but with $\newprodinline$ not
equal to the Arens products.
\end{proposition}
\begin{proof}
Let our nets be indexed by the directed set $J$, and let $\mc V$ be an
ultrafilter refining the order filter on $J$.  Let $\mc W$ be some ultrafilter
such that there exists a bounded map $L:\mc A''\rightarrow(\mc A)_{\mc W}$
satisfying the conclusions of Theorem~\ref{arens_case}.  Let $\Phi\in\mc A''$ be the
weak$^*$-limit of $(a_\alpha b_\alpha)$ along $\mc V$, so that $\Phi\not=0$.

Let $\Psi_1,\Psi_2\in\mc A''\setminus\mc A$ and $M_1,M_2\in\mc A'''$
be such that $\ip{M_1}{\kappa_{\mc A}(a)} =
\ip{M_2}{\kappa_{\mc A}(a)} = 0$ for $a\in\mc A$, and
\[ \ip{M_1}{\Psi_1} = \ip{M_2}{\Psi_2} = 1, \quad
\ip{M_2}{\Psi_1} = \ip{M_1}{\Psi_2} = 0. \]
Recall that we define the ultrafilter $\mc V\times\mc W$ on $J\times I$
by, for $K\subseteq J \times I$, setting $K\in\mc V\times\mc W$ if and only if
\[ \big\{ i\in I : \{ \alpha\in J : (\alpha,i)\in K \} \in \mc V \big\} \in \mc W. \]
Then, for a family $(x_{\alpha,i})_{\alpha\in J, i\in I}$ in a Hausdorff
space $X$, we have that
\[ \lim_{i\rightarrow\mc W} \lim_{\alpha\rightarrow\mc V} x_{\alpha,i} =
\lim_{(\alpha,i)\rightarrow\mc V\times\mc W} x_{\alpha,i}, \]
whenever the limits exist.
Let $\mc U = \mc V\times\mc W$, and define
$K:\mc A''\rightarrow(\mc A)_{\mc U}$ as follows.  For $\Lambda\in\mc A''$,
let $L(\Lambda) = (c_i)_{i\in I}$, and define
\[ K(\Lambda) = \big(c_i + \ip{M_1}{\Lambda}a_\alpha + \ip{M_2}{\Lambda}b_\alpha
\big)_{(\alpha,i)\in J\times I}. \]
Obviously $K$ is linear and bounded.  For $a\in\mc A$, by the choice of
$M_1$ and $M_2$, we have that $K(\kappa_{\mc A}(a)) = (a)$, so that
$K\circ\kappa_{\mc A}$ is the canonical map $\mc A\rightarrow(\mc A)_{\mc U}$.

For $\mu\in\mc A'$ and $\Lambda\in\mc A''$, let $L(\Lambda)=(c_i)$, so we have that
\begin{align*} \ip{\sigma(K(\Lambda))}{\mu} =
\lim_{(\alpha,i)\rightarrow\mc U} \ip{\mu}{c_i + \ip{M_1}{\Lambda}a_\alpha
   + \ip{M_2}{\Lambda}b_\alpha}
= \lim_{i\rightarrow\mc W} \ip{\mu}{c_i} = \ip{\Lambda}{\mu}, \end{align*}
as $(a_\alpha)$ and $(b_\alpha)$ are weakly-null.

Finally, let $L(\Psi_1)=(c_i)$ and $L(\Psi_2)=(d_i)$, so that
\[ K(\Psi_1) = (c_i + a_\alpha), \quad K(\Psi_2) = (d_i+b_\alpha). \]
Thus we see that for $\mu\in\mc A'$,
\begin{align*} \ip{\Psi_1 \newprod \Psi_2}{\mu} &=
\lim_{(\alpha,i)\rightarrow\mc U} \ip{\mu}{(c_i + a_\alpha)(d_i + b_\alpha)}
= \lim_{(\alpha,i)\rightarrow\mc U} \ip{\mu}{c_id_i + a_\alpha d_i + c_i b_\alpha + a_\alpha b_\alpha} \\
&= \ip{\Psi_1\aone\Psi_2}{\mu} + \ip{\Phi}{\mu} +
\lim_{i\rightarrow\mc W}\Big( \lim_{\alpha\rightarrow\mc V} \ip{d_i\cdot\mu}{a_\alpha}
   + \ip{\mu\cdot c_i}{b_\alpha} \Big) \\
&= \ip{\Psi_1\aone\Psi_2 + \Phi}{\mu},
\end{align*}
as $(a_\alpha)$ and $(b_\alpha)$ are weakly-null.
Hence $\Psi_1 \newprod \Psi_2 \not= \Psi_1\aone\Psi_2$ as $\Phi\not=0$.
\end{proof}

Let $\mc A = \mc K(\ell^p)$, the compact operators of $\ell^p$, for
$1<p<\infty$.  Then, as detailed in \cite[Section~1.4]{Palmer1} for example,
the dual of $\mc A$ is the nuclear operators on $\ell^p$, and
$\mc A'' = \mc B(\ell^p)$, with $\aone=\atwo$ agreeing with the usual
product.  Let $(\delta_n)$ be the standard unit vector basis of
$\ell^p$, and for each $n$, let $a_n$ be the rank-one operator
which sends $\delta_n$ to $\delta_1$, and kills $\delta_k$ otherwise,
and let $b_n$ be the rank-one operator which sends $\delta_1$ to $\delta_n$,
and kills $\delta_k$ otherwise.  Then $(a_n)$ and $(b_n)$ are both
weakly-null sequences, but $a_nb_n$ is the projection onto the first
co-ordinate, for each $n$.  Hence the above proposition applies
in this case, and in particular it applies to the C$^*$-algebra
$\mc K(\ell^2)$.

\begin{definition}
Let $E$ and $F$ be Banach spaces, and let $T:E\rightarrow F$ be a bounded
linear map.  Then $T$ is \emph{totally completely continuous} if, whenever
$(x_\alpha)$ is a weakly-null net in $E$, then $\lim_\alpha \| T(x_\alpha) \|=0$.
\end{definition}

Recall that if we only use sequences in the above definition, we get the
usual notion of $T$ being \emph{completely continuous}.  Our work below will
show that being totally completely continuous is strictly stronger than
being completely continuous.  We view being totally completely continuous
as a property close to being compact.

\begin{proposition}
Let $\mc A$ be a Banach algebra, and let $\mu\in\wap(\mc A')$ be such
that the map $\mc A\rightarrow\mc A'; a\mapsto a\cdot\mu$ is not
totally completely continuous.  Then there exist weakly-null nets
$(a_\alpha)$ and $(b_\alpha)$, on the same index set,
such that $(b_\alpha a_\alpha )$ is not weakly-null.
\end{proposition}
\begin{proof}
Let $(a_i)$ be some weakly-null net in $\mc A$ such that for
some $\delta>0$, we have that $\|a_i\cdot\mu\|\geq\delta$ for all $i$.
Let $\Lambda$ be the collection of finite-dimensional subspaces of
$\mc A'$, partially ordered by reverse inclusion.  For each $M\in\Lambda$,
if we can find some $i$ and some $b\in\mc A$ with $\|b\|=1$, and with
\[ |\ip{a_i\cdot\mu}{b}| \geq \delta/3, \qquad
\ip{\lambda}{b} = 0, \quad
|\ip{\lambda}{a_i}|\leq (\dim M)^{-1} \|\lambda\|
\qquad (\lambda\in M), \]
then we are done.

As $(a_i)$ is weakly-null, we can ensure the final condition by
simply ensuring that $i\geq i_0$, for some $i_0$ depending upon $M$.
For $c\in\mc A$, we have that
\[ \lim_i \ip{a_i\cdot\mu}{c} = \lim_i \ip{\mu\cdot c}{a_i} = 0, \]
so we see that $(a_i\cdot\mu)$ is weak$^*$-null.

For fixed $i$, notice that the Hahn-Banach theorem shows that
\begin{align*} \sup\big\{ |\ip{a_i\cdot\mu}{b}| &: b\in\mc A, \|b\|=1,
   \ip{\lambda}{b}=0 \ (\lambda\in M) \big\} \\
&= d(a_i\cdot\mu,M) := \inf\big\{ \|a_i\cdot\mu-\lambda\| : \lambda\in M \big\}.
\end{align*}
For each $i$, let $\lambda_i\in M$ be such that $\|a_i\cdot\mu-\lambda_i\|
=d(a_i\cdot\mu,M)$.  Suppose that $d(a_i\cdot\mu,M)\leq\epsilon$ for each
$i\geq i_0$.  Then $(\lambda_i)$ is a bounded net, and so, by
passing to a subnet if necessary, we may suppose that $\lambda_i\rightarrow\lambda
\in M$, in norm.  For $c\in\mc A$, we see that
\[ |\ip{\lambda}{c}| = \lim_i |\ip{a_i\cdot\mu-\lambda}{c}|
= \lim_i |\ip{a_i\cdot\mu-\lambda_i}{c}| \leq \|c\| \epsilon, \]
so that $\|\lambda\|\leq\epsilon$.  Hence 
\[ \delta \leq \liminf_i \|a_i\cdot\mu\|
\leq \liminf_i \|a_i\cdot\mu-\lambda\| + \|\lambda\|
\leq \epsilon + \liminf_i \|a_i\cdot\mu-\lambda_i\| \leq 2\epsilon, \]
showing that $\epsilon\geq \delta/2$, as required.
\end{proof}

For further about the following definition, see \cite[Section~5]{Cas}.

\begin{definition}
A Banach space $E$ has the $\pi$-property if there exists a
bounded net $(T_\alpha)$ of finite-rank projections such that
$\|T_\alpha(x)-x\|\rightarrow0$ for $x\in E$.
\end{definition}

Notice that if $E$ has a basis, then it has the $\pi$-property.
It is easy to show that $L^1$ spaces (and hence $M(X)$ spaces)
have the $\pi$-property, so the duals of commutative C$^*$-algebras
have the $\pi$-property.  However, the $\pi$-property is stronger
than the more usual (bounded) approximation property, and it is
known (see \cite[Remark~6.1.9]{RundeBook}) that a von Neumann algebra
$\mc M$ has the approximation property if and only if it is nuclear,
which is if and only if $\mc M$ is a finite sum of algebras of the
form $C(X)\otimes\mathbb M_n$.  As the approximation property
passes to preduals, there are plenty of von Neumann algebras $\mc M$
such that $\mc M'$ cannot have the $\pi$-property.

\begin{proposition}
Let $\mc A$ be a Banach algebra such that $\mc A'$ has the $\pi$-property.
Let $\mu\in\wap(\mc A')\setminus\ap(\mc A')$.  Then the map
$\mc A\rightarrow\mc A'; a\mapsto a\cdot\mu$ is not
totally completely continuous.
\end{proposition}
\begin{proof}
Let $(T_\alpha)$ be a net of finite-rank projections on $\mc A'$ such that
$\|T_\alpha(\lambda)-\lambda\|\rightarrow0$ for each $\lambda\in\mc A'$.
Let $M = \sup_\alpha \|T_\alpha\|$.
Suppose, towards a contradiction, that for each $\delta>0$ there exists
some $\alpha$ such that for each $\Phi\in\mc A''$ with $T_\alpha'(\Phi)=0$,
we have that $\|\Phi\cdot\mu\|\leq\delta\|\Phi\|$.

For each $\delta$, choose such an $\alpha$.  For $\Phi\in\mc A''$ let
$\Psi = \Phi - T_\alpha'(\Phi)$, so that $T_\alpha'(\Psi)=0$, and so
$\|\Phi\cdot\mu - T_\alpha'(\Phi)\cdot\mu\| = \|\Psi\cdot\mu\| \leq
\delta\|\Psi\| \leq \delta (1+M) \|\Phi\|$.  As $\delta>0$ was arbitrary,
we see that the map $\mc A''\rightarrow\mc A'; \Phi\mapsto \Phi\cdot\mu$
can be uniformly approximated by finite-rank operators, and hence is compact,
contradicting the fact that $\mu\not\in\ap(\mc A')$.

Thus, for some $\delta>0$, for each $\alpha$ there exists $\Phi_\alpha\in
\mc A''$ with $\|\Phi_\alpha\|=1$, $T_\alpha'(\Phi_\alpha)=0$, and
$\|\Phi_\alpha\cdot\mu\|\geq\delta$.  For each $\alpha$, let $a_\alpha
\in\mc A$ be such that $\|a_\alpha\|\leq2$, $a_\alpha$ agrees with $\Phi_\alpha$
on the image of $T_\alpha$, a finite-dimensional subspace of $\mc A'$, and
$\|a_\alpha\cdot\mu\|\geq\delta/2$.  We can ensure the final condition
as $a_\alpha\cdot\mu \rightarrow \Phi_\alpha$ weakly, as $\mu\in\wap(\mc A')$.
Then, for $\lambda\in\mc A'$, as $T_\alpha(\lambda)\rightarrow\lambda$ in
norm, we see that
\[ \lim_\alpha \ip{\lambda}{a_\alpha}
= \lim_\alpha \ip{T_\alpha(\lambda)}{a_\alpha}
= \lim_\alpha \ip{\Phi_\alpha}{T_\alpha(\lambda)} = 0, \]
so $(a_\alpha)$ is weakly-null.
\end{proof}

We can immediately draw conclusions about commutative C$^*$-algebras.

\begin{theorem}
Let $\mc A = C_0(X)$ be a commutative C$^*$-algebra.  The following
are equivalent:
\begin{enumerate}
\item $X$ is scattered;
\item for any $K:\mc A''\rightarrow(\mc A)_{\mc U}$ satisfying
our standing assumptions, we have that $\newprodinline = \aone$.
\end{enumerate}
\end{theorem}

Asking for $K$ to be an isometry is a strong condition.  There is by now
a large selection of results in the theory of C$^*$-algebras which weaken
the requirement of ``isometry'' to ``completely-bounded'' (see \cite{Paul},
for example).  Given a C$^*$-algebra $\mc A$, there is a canonical way
to turn the matrix algebra $\mathbb M_n(\mc A)$ into a C$^*$-algebra.
Given a linear map $T:\mc A\rightarrow\mc B$ between two C$^*$-algebras,
we let $T$ act pointwise as a map $(T)_n:\mathbb M_n(\mc A)\rightarrow\mathbb M_n(\mc B)$.
Then $T$ is \emph{completely-bounded} if and only if $\sup_n \|(T)_n\| <\infty$.
It is simple to show the map $K$ which we construct in
Proposition~\ref{when_K_not_good} is completely-bounded: the only
non-trivial thing to check is the well-known fact that bounded linear functionals
are automatically completely-bounded.  It would be interesting to know if
there is any reasonable weakening of the ``isometry'' condition on $K$ which
still ensures that $\newprodinline=\aone$ for C$^*$-algebras.

\section{Acknowledgments}

Some of the results in this paper are from the author's
PhD thesis \cite{DawsThesis} completed at the University of
Leeds under the financial support the the EPSRC.  The author
wishes to thank his PhD supervisors Garth Dales and Charles Read.

\vspace{5ex}

\noindent\emph{Author's Address:}
\parbox[t]{3in}{St. John's College,\\
Oxford,\\
OX1 3JP.}

\bigskip\noindent\emph{Email:} \texttt{matt.daws@cantab.net}


\begin{thebibliography}{99}
\normalsize
\newcommand{\bibbook}[3]{\textsc{#1}, \emph{#2}, (#3).}
\newcommand{\bibpaper}[6]{\textsc{#1}, `#2', \emph{#3} #4 (#5) #6.}
\newcommand{\bibpreprint}[2]{\textsc{#1}, `#2', preprint.}

\bibitem{AK} \bibbook{F. Albiac, N.\,J. Kalton}
   {Topics in Banach space theory}
   {Springer, New York, 2006}

\bibitem{Arens} \bibpaper{R. Arens}
   {The adjoint of a bilinear operation}
   {Proc. Amer. Math. Soc.}{2}{1951}{839--848}

\bibitem{Beh} \bibpaper{E. Behrends}
   {On the principle of local reflexivity}
   {Studia Math.}{100}{1991}{109--128}

\bibitem{Cas} \textsc{P.\,G. Casazza},
   `Approximation properties', \emph{Handbook of the geometry of Banach spaces, Vol.\ I},
   271--316, (North-Holland, Amsterdam, 2001)

\bibitem{Dales} \bibbook{H.\,G. Dales}
   {Banach algebras and automatic continuity}
   {Clarendon Press, Oxford, 2000}

\bibitem{DawsDBA} \bibpaper{M. Daws}
   {Dual Banach algebras: representations and injectivity}
   {Studia Math.}{178}{2007}{231--275}

\bibitem{DawsThesis} \textsc{M. Daws}, \emph{Banach algebras of operators},
   PhD. thesis, University of Leeds, 2005.

\bibitem{DF} \bibbook{A. Defant, K. Floret}
   {Tensor norms and operator ideals}
   {North-Holland Publishing Co., Amsterdam, 1993}

\bibitem{Diestel} \bibbook{J. Diestel}
   {Sequences and series and Banach spaces}
   {Springer-Verlag, New York, 1984}

\bibitem{DU} \bibbook{J. Diestel, J.J. Uhl, Jr.}
   {Vector measures}{American Mathematical Society, Providence, 1977}

\bibitem{GH} \textsc{L. Ge, D. Hadwin},
   `Ultraproducts of $C\sp *$-algebras', in
   \emph{Recent advances in operator theory and related topics}, 305--326,
   Oper. Theory Adv. Appl., 127, (Birkhäuser, Basel, 2001).

\bibitem{GI} \bibpaper{G. Godefroy, B. Iochum}
   {Arens-regularity of Banach algebras and the
   geometry of Banach spaces}
   {J. Funct. Anal.}{80}{1988}{47--59}

\bibitem{Hein} \bibpaper{S. Heinrich}{Ultraproducts in Banach space theory}
	{J. reine angew. Math.}{313}{1980}{72--104}

\bibitem{IL} \bibpaper{B. Iochum, G. Loupias}
   {Arens regularity and local reflexivity
   principle for Banach algebras}
   {Math. Ann.}{284}{1989}{23--40}

\bibitem{Lau} \bibpaper{A.\,T.-M. Lau}
   {Uniformly Continuous Functionals on the Fourier Algebra of Any Locally Compact Group}
   {Trans. Amer. Math. Soc.}{251}{1979}{39--59}

\bibitem{LL2} \bibpaper{A.\,T.-M. Lau, R.\,J. Loy}
   {Weak amenability of Banach algebras on locally compact groups}
   {J. Funct. Anal.}{145}{1997}{175--204}

\bibitem{Palmer1} \bibbook{T.\,W.~Palmer}
   {Banach algebras and the general theory of $\sp *$-algebras. Vol. 1.}
   {Cambridge University Press, Cambridge, 1994}

\bibitem{Paul} \bibbook{V. Paulsen, Vern}
   {Completely bounded maps and operator algebras}
   {Cambridge University Press, Cambridge, 2002}

\bibitem{Quigg} \bibpaper{J.\,C. Quigg}
   {Approximately periodic functionals on $C\sp *$-algebras and von Neumann algebras}
   {Canad. J. Math.}{37}{1985}{769--784}

\bibitem{RundeBook} \bibbook{V. Runde}
   {Lectures on amenability}
   {Springer-Verlag, Berlin, 2002}

\bibitem{Runde} \bibpaper{V. Runde}
   {Amenability for dual Banach algebras}
   {Studia Math.}{148}{2001}{47--66}

\bibitem{Ryan} \bibbook{R. Ryan}
   {Introduction to Tensor Products of Banach Spaces}
   {Springer-Verlag, London, 2002}

\bibitem{Tak} \bibbook{M. Takesaki}
   {Theory of Operator Algebras I}
   {Springer-Verlag, New York, 1979}


\end{thebibliography}
\end{document}